\documentclass[11pt]{amsart}
\usepackage{amsmath,amsfonts,amsthm,amssymb,amscd}
\def\classification#1{\def\@class{#1}}
\classification{\null}

\DeclareFontFamily{OT1}{rsfs}{}
\DeclareFontShape{OT1}{rsfs}{n}{it}{<-> rsfs10}{}
\DeclareMathAlphabet{\mathscr}{OT1}{rsfs}{n}{it}

\DeclareMathOperator{\diam}{diam}
\DeclareMathOperator{\SL}{SL}

\DeclareMathOperator{\PSL}{PSL}
\DeclareMathOperator{\SU}{SU}

\DeclareMathOperator{\mix}{mix}

\DeclareMathOperator{\Tr}{Tr}
\DeclareMathOperator{\tr}{tr}

\newtheorem{prop}{Proposition}[section]
\newtheorem{thm}[prop]{Theorem}
\newtheorem*{main}{Main Theorem}
\newtheorem*{kprop}{Key Proposition}
\newtheorem*{conj}{Conjecture}
\newtheorem{cor}[prop]{Corollary}
\newtheorem{lem}[prop]{Lemma}
\newtheorem{defn}{Definition}

\numberwithin{equation}{section}
\title{Growth and generation in $\SL_2(\mathbb{Z}/p \mathbb{Z})$}
\author{H. A. Helfgott}
\address{H. A. Helfgott, Mathematics Department, University of Bristol, University Walk, Bristol BS8 1TW, United Kingdom}
\email{h.andres.helfgott@bristol.ac.uk}
\subjclass[2000]{05C25, 20G40, 20D60, 11B75}
\keywords{Cayley graphs, finite groups, generation, diameter}
\thanks{The author was supported by a fellowship from the Centre de
Recherches Math\'ematiques at Montr\'eal. Travel was partially 
funded by the Clay Mathematics Institute.} 
\begin{document}
\begin{abstract}
We show that every subset of $\SL_2(\mathbb{Z}/p\mathbb{Z})$ grows 
rapidly when it acts on itself by the group operation.
 It follows readily that, for
every set of generators $A$ of $\SL_2(\mathbb{Z}/p\mathbb{Z})$, every
element of $\SL_2(\mathbb{Z}/p\mathbb{Z})$ can be expressed as a product of
at most $O((\log p)^c)$ elements of $A \cup A^{-1}$, where
$c$ and the implied constant are absolute.
\end{abstract}
\maketitle
\section{Introduction}
\subsection{Background}
Let $G$ be a finite group. Let $A\subset G$ be a set of generators of $G$.
By definition, every $g\in G$ can be expressed as a product of elements
of $A\cup A^{-1}$. We would like to know the length of the
longest product that might be needed; in other words, we wish to bound
from above the diameter $\diam(\Gamma(G,A))$ of the Cayley graph of $G$ with 
respect to $A$. (The {\em Cayley graph} $\Gamma(G,A)$ is the graph
$(V,E)$ with vertex set $V = G$ and edge set $E = \{(a g, g) :
g\in G, a\in A\}$. The {\em diameter} of a graph $X = (V,E)$
is $\max_{v_1,v_2\in V} d(v_1,v_2)$, where $d(v_1,v_2)$ is the length of
the shortest path between $v_1$ and $v_2$ in $X$.)

If $G$ is abelian, the diameter can be very large: if $G$ is cyclic of
order $2 n + 1$, and $g$ is any generator of $G$, then $g^n$ cannot be
expressed as a product of length less than $n$ on the elements of
$\{g, g^{-1}\}$. However, if $G$ is non-abelian and simple, the diameter
is believed to be quite small:

\begin{conj}[Babai, \cite{BS}] For every non-abelian finite simple group $G$,
\begin{equation}\label{eq:udo}
\diam(\Gamma(G,A)) \ll (\log |G|)^c,\end{equation}
where $c$ is some absolute constant and $|G|$ is the number of elements of $G$.
\end{conj}

This conjecture is far from being proved. Even for the basic cases, viz.,
$G = A_n$ and $G = \PSL_2(\mathbb{Z}/p\mathbb{Z})$, the conjecture 
has remained open
until now; these two choices of $G$ seem to present already many of the main 
difficulties of the general case.

Work on both kinds of groups long predates the general conjecture
in \cite{BS}.
Let us focus\footnote{
While $\SL_2(\mathbb{Z}/p\mathbb{Z})$ is not simple, the statement
(\ref{eq:udo}) for $\SL_2(\mathbb{Z}/p\mathbb{Z})$ is trivially
equivalent to (\ref{eq:udo}) for $\PSL_2(\mathbb{Z}/p\mathbb{Z})$,
and treating the former group is both slightly more conventional
and notationally simpler.}
on $G = \SL_2(\mathbb{Z}/p\mathbb{Z})$.
 There are some classical
results for certain specific generators. Let
\begin{equation}\label{eq:knop}
A = \left\{ \left(\begin{array}{cc} 1 &1\\0 &1\end{array}\right),
\left(\begin{array}{cc} 1 &0\\1 &1\end{array}\right)
\right\} .\end{equation}
Selberg's spectral-gap theorem
for $\SL_2(\mathbb{Z})\backslash \mathbb{H}$ (\cite{Se}) implies that
$\{\Gamma(\SL_2(\mathbb{Z}/p \mathbb{Z}), A)\}_{p\geq 5}$ 
is a family
of expander graphs (vd., e.g., \cite{Lu}, Thm.\ 4.4.2, (i)). 
It follows easily that
\[\diam(\Gamma(\SL_2(\mathbb{Z}/p \mathbb{Z}),A))\ll \log p .\] 
Unfortunately, this argument works only for a few
 other choices of $A$. For example, no good bounds were known up to now
for $\diam(\Gamma(\SL_2(\mathbb{Z}/p\mathbb{Z}),A))$ with, say,
\begin{equation}\label{eq:ujo}
A = \left\{ \left(\begin{array}{cc} 1 &3\\0 &1\end{array}\right),
\left(\begin{array}{cc} 1 &0\\3 &1\end{array}\right)
\right\} ,\end{equation}
let alone for general $A$, uniformly on $A$ or not.
\subsection{Results}\label{subs:resu}
We prove the conjecture 
for $G = \SL_2(\mathbb{Z}/p \mathbb{Z})$.
\begin{main}
Let $p$ be a prime. Let $A$ be a set of generators of
$G = \SL_2(\mathbb{Z}/p \mathbb{Z})$. Then the Cayley graph
$\Gamma(G,A)$ has diameter $O((\log p)^c)$, where
$c$ and the implied constant are absolute.
\end{main} 
The theorem is a direct consequence of the following statement.
\begin{kprop}
Let $p$ be a prime. Let $A$ be a subset of $\SL_2(\mathbb{Z}/p \mathbb{Z})$
not contained in any proper
subgroup.
\begin{enumerate}
\item\label{it:wark} Assume that $|A| < p^{3 - \delta}$ for some
fixed $\delta>0$.
Then
\begin{equation}\label{eq:solt}
|A\cdot A\cdot A| > c |A|^{1 +\epsilon},\end{equation}
where $c>0$ and $\epsilon>0$ depend only on $\delta$.
\item\label{it:agust}
Assume that $|A| > p^{\delta}$ for some fixed $\delta>0$. Then there 
is an integer $k>0$, depending only on $\delta$, such that
every element of $\SL_2(\mathbb{Z}/p \mathbb{Z})$ can be expressed
as a product of at most $k$ elements of $A \cup A^{-1}$.
\end{enumerate}
\end{kprop}
The crucial fact here is that the constants $c$, $\epsilon$ and $k$ do 
not depend on $p$ or on $A$.

It follows immediately from the main theorem (via \cite{DSC}, \S 2, Lem.\ 2, 
\S 3, Cor.\ 3.1, and \S 3, Cor.\ 3.2) that
the {\em mixing time} of $\Gamma(\SL_2(\mathbb{Z}/p \mathbb{Z}),A)$ is 
$O(|A| (\log p)^{2 c + 1})$,
where $c$ and the implied constant are absolute, and $c$ is as in the
main theorem. (The {\em mixing time} is
the least $t$ for which a lazy 
random walk of length $t$ starting at the origin of
the Cayley graph has a distribution of destinations close to the uniform
distribution in the $\ell_1$ norm; vd.\ \S \ref{sec:conc})

If $A$ equals the projection of a fixed set of generators of a free
group in $\SL_2(\mathbb{Z})$ (take, e.g., $A$ as in (\ref{eq:knop}) or 
(\ref{eq:ujo})) it follows by a simple argument that $A$ must grow rapidly
at first when multiplied by itself. In such a situation, we obtain a bound
of 
\[\diam(\Gamma(\SL_2(\mathbb{Z}/p \mathbb{Z}),A)) \ll \log p,\]
where the implied constant depends on the elements of $\SL_2(\mathbb{Z})$
of which $A$ is a projection. 
For (\ref{eq:ujo}) and most other examples, this bound is new; for
$A$ as in (\ref{eq:knop}), it is, of course, known, and the novelty lies in 
the proof\footnote{What is given here is not, however, the first 
elementary proof for the choice of $A$ in (\ref{eq:knop}); see
\cite{SX}. The proof in \cite{SX} works for all projections of sets generating
finite-index subgroups of $\SL_2(\mathbb{Z})$. Gamburd \cite{Gambo}
succeeded in extending the method to projections of sets generating subgroups
of $\SL_2(\mathbb{Z})$ whose limit sets have
 Hausdorff dimension greater than $5/6$.}.

If $A$ is a random pair of generators, then,
with probability tending to $1$ as $p\to \infty$, the graph
$\Gamma(\SL_2(\mathbb{Z}/p \mathbb{Z}),A)$ does not have small loops (see \S \ref{sec:conc}). 
It then follows from the key
proposition that $\diam(\Gamma(\SL_2(\mathbb{Z}/p \mathbb{Z}),A)) \ll \log p$,
as ventured by Lubotzky (\cite{Lu}, Prob.\ 10.3.3). The implied constant
is absolute.
\subsection{Techniques}
The tools used are almost exclusively additive-combinatorial. Fourier
analysis over finite fields and Ruzsa distances are used repeatedly.
Both Gowers's effective version of the Balog-Szemer\'edi theorem
(\cite{Gow}) and the sum-product estimates in \cite{BKT} and
\cite{Ko} play crucial roles. It is only through \cite{Ko} that 
arithmetic strictly speaking plays a role, viz., in the guise
of an estimate proved in \cite{HBK} with techniques derived from
Stepanov's elementary proof of the Weil bounds. The Weil bounds themselves
are not used, and even the use of \cite{Ko} becomes unnecessary when
auxiliary results suffice to ensure the growth of $A$ small (namely, in
the cases of fixed or random generators).

Estimates on growth in $\mathbb{Z}/p \mathbb{Z}$ 
will be proved in \S \ref{sec:exf},
and part
(\ref{it:wark}) of the
 key proposition will be reduced thereto in \S \ref{sec:gerg}. Given
part (\ref{it:wark}), it suffices to prove (\ref{it:agust}) for very
large $A$ -- and this is a relatively simple task (\S \ref{sec:whg}), 
yielding to the use of growth estimates coming from Fourier analysis.

\subsection{Work to do}
A natural next step would be to generalise the main results to the group
$\SL_2(\mathbb{F}_{p^{\alpha}})$, $\alpha>1$. At first sight, this does not
seem too hard; however, there seem to be actual difficulties in making
the result uniform on $\alpha$. 

A generalisation to $\SL_n(\mathbb{Z}/p\mathbb{Z})$ for $n\geq 3$ is likely
to require a great deal of original work. The arguments in 
\S \ref{subs:noeth}-\ref{subs:sitra} should carry over, but those in
\S \ref{sec:exf} and \S \ref{subs:grosma} do not. It is possible
that the basic approach in \S \ref{subs:noeth}-\ref{subs:sitra} will
eventually prove itself valid for all simple\footnote{
The diameter of a Cayley graph $\Gamma(G,A)$ of a solvable linear
algebraic group $G$ can be large: for example, $G$ could be generated by
the set $A$ of all elements of $G$ all of whose eigenvalues lie in
$B$, where 
$B\subset (\mathbb{F}_{p^{\alpha}})^*$ is a set that
 grows very slowly when multiplied by itself. By the Lie-Kolchin
theorem, the eigenvalues of $A\cdot A\dotsb A$ will lie in
$B\cdot B\dotsb B$, which, by assumption, is only slightly larger than $B$.
(See also \cite{ET}.)
It is unclear whether the present paper's approach will be directly
applicable to groups that are neither solvable nor simple 
(nor almost simple).}  groups of Lie type, but it
is too soon to tell whether something will be found to replace
\S \ref{sec:exf} and \S \ref{subs:grosma} in a general context.


No attempt has been made to optimize -- or compute -- the constant $c$
in the main theorem, though, like the implied constant, it is effective
and can be made explicit. 
Actual numerical constants will sometimes be used in the argument for the
sake of notational clarity.

\subsection{Further remarks}
There is a rich literature on the growth of sets in linear algebraic
groups over fields of characteristic zero: see, most recently, \cite{EMO}. 
In such a situation, one has access to topological arguments without 
clear analogues in $\mathbb{Z}/p \mathbb{Z}$. It is possible, nevertheless,
to adapt the vocabulary of growth on infinite groups to the finite case. For 
example, one can say the key proposition implies immediately
that $A$ does not have
{\em moderate growth} (\cite{DSC2}). 

The problem of bounding the diameter of
$\Gamma(\SL_2(\mathbb{Z}/p^k \mathbb{Z}),A)$ for $p$ fixed and $k$ variable
 is fundamentally different from that of bounding the diameter 
of $\Gamma(\SL_2(\mathbb{Z}/p \mathbb{Z}), A)$ for $p$ variable.
From a $p$-adic perspective, the problem for 
$\SL_2(\mathbb{Z}/p^k \mathbb{Z})$ is analogous to that for
$\SU(2)$, which was treated by Solovay and Kitaev \cite{NC}.
Dinai \cite{Di} has succeeded in giving a polylogarithmic bound for
$\diam(\Gamma(\SL_2(\mathbb{Z}/p^k \mathbb{Z}),A))$, $p$ fixed,
in part by adapting Solovay and Kitaev's procedure.

Consider the family 
$\mathscr{F} = \{\Gamma(\SL_2(\mathbb{Z}/p \mathbb{Z}), A)\}_{
p,A}$, where both $p$ and $A$ vary: $p$ ranges across the primes and $A$
ranges across all sets that generate $\SL_2(\mathbb{Z}/ p \mathbb{Z})$.
If we could prove that $\mathscr{F}$ is an expander family, we would obtain
the main theorem with the constant $c$ set to $1$. We are still
far from proving that $\mathscr{F}$ is an expander family, and we will not,
of course, assume such a hypothesis;
 rather, we will obtain a weaker statement as an 
immediate consequence of the main theorem (Cor.\ \ref{cor:cory}).
It seems unjustified for now to hope for
a purely combinatorial proof that a family of Cayley graphs $\{\Gamma(G,A)\}$
where both $G$ and $A$ vary quite freely is an
expander family: we would need, not estimates on the growth of a set
$A$ when added to or multiplied by itself, but, instead, estimates on
the growth of a set $A$ under the action of addition or multiplication
by a small, fixed set $S$, or under the action of a small set of operations.
(Here ``small'' means ``of cardinality less than a constant''.)
Such estimates are outside of the reach of the already remarkably
strong sum-product techniques of \cite{BKT} and \cite{Ko}.
\subsection{Acknowledgments}
I would like to thank A. Venkatesh for
having first called the problem to my attention and for
shedding light spontaneously. His Clay Mathematics
Institute grant paid for a trip during which the present subject
and many other interesting things were discussed.
I was otherwise funded by the Centre de Recherches Math\'ematiques
and the Institut de Sciences Math\'ematiques (Montr\'eal).

Thanks are also due to N.\ Anantharaman,
E.\ Breuillard, O.\ Dinai, U.\ Hadad, C.\ Hall and
G.\ Harcos, for their careful reading and several helpful
comments, to A.\ Gamburd, A.\ Lubotzky and I.\ Pak, for their
instructive remarks and references,
and to A.\ Granville, for his encouragement and advice, and for
access to an unpublished set of lecture notes.

\section{Background and preliminaries}
\subsection{General notation}
As is customary, we denote by $\mathbb{F}_{p^{\alpha}}$ the finite field
of order $p^{\alpha}$. We write $|f|_r$ for the $L_r$--norm of a function
$f$. Given a set $A$, we denote its cardinality by $|A|$, and its
characteristic function by $A$ itself. Thus, $|A| = |A|_1$. By $A+B$
(resp. $A\cdot B$), we shall always mean
$\{x+y : x\in A, y\in B\}$ (resp. $\{x\cdot y: x\in A, y\in B\}$), 
or the characteristic
function thereof; cf. $(A\ast B) (x) = 
|\{(y,z) \in A\times B : y+z = x\}|$. By $A+\xi$ and $\xi\cdot A$ we 
mean $\{x+\xi : x\in A\}$ and $\{\xi \cdot x: x\in A\}$, respectively.

For us, $A^r$ means 
$\{x^r : x\in A\}$; in general, if $f$ is a function on $A$,
we take $f(A)$ to mean $\{f(x) : x\in A\}$.
Given a positive integer $r$ and a
 subset $A$ of a group $G$, we define
$A_r$ to be the set of all products of at most $r$ elements of $A \cup A^{-1}$:
\[A_r = \{g_1 \cdot g_2 \dotsb g_r : g_i\in A \cup A^{-1} \cup \{1 \}\} .\]
Finally, we write $\langle A\rangle$ for the group generated by $A$.
\subsection{Fourier analysis over $\mathbb{Z}/p\mathbb{Z}$}
We will review some basic facts, in part to fix our normalizations.
The Fourier transform $\widehat{f}$ of a function $f:\mathbb{Z}/p\mathbb{Z}\to
\mathbb{C}$ is given by
\[\widehat{f}(y) = \sum_{x\in \mathbb{Z}/p\mathbb{Z}} f(x) e^{-2 \pi i x y/p} .\]
The Fourier transform is an isometry:
\[\sum_{x\in \mathbb{Z}/p\mathbb{Z}} |\widehat{f}(x)|^2 = p \cdot
\sum_{x\in \mathbb{Z}/p\mathbb{Z}} |f(x)|^2 .\]
For any $f,g: \mathbb{Z}/p\mathbb{Z} \to \mathbb{C}$, we have 
$\widehat{f\ast g} = \widehat{f} \cdot \widehat{g}$. If $A,B \subset
\mathbb{Z}/p\mathbb{Z}$, then $|A\ast B|_1 = |A| |B|$.
\subsection{Additive combinatorics, abelian and non-abelian}
Some basic concepts and proofs of additive combinatorics transfer effortlessly
to the non-abelian case; some do not. In the following, $G$ need not be
an abelian group, except, of course, when it is explicitly said to be one.
\begin{defn}
Let $A$ and $B$ be finite
subsets of a group $G$. We define the {\em Ruzsa distance}
\[d(A,B) = \log\left(\frac{|A B^{-1}|}{\sqrt{|A| |B|}}\right) .\]
\end{defn}
If $G$ is an abelian group whose operation is written additively, we denote
the Rusza distance by $d_+(A,B)$.

The Ruzsa distance, while not truly a distance function ($d(A,A)\ne 0$ in
general), does satisfy the triangle inequality.
\begin{lem}\label{lem:ineq}
Let $A$, $B$ and $C$ be finite subsets of a group $G$.
Then
\begin{equation}\label{eq:trian}
d(A,C) \leq d(A,B) + d(B,C) .\end{equation}
\end{lem}
\begin{proof}[Proof (Ruzsa)]
It is enough to prove that
\begin{equation}\label{eq:escolt}
|A C^{-1}| |B| \leq |A B^{-1}| |B C^{-1}| .\end{equation}
We will do as much by constructing an injection
$\iota:A C^{-1} \times B \hookrightarrow A B^{-1} \times
 B C^{-1}$. For every 
$d\in A C^{-1}$, choose once and for all a pair
$(a_d,c_d)\in A\times C$ such that
$d = a_d c_d^{-1}$. Define 
$\iota(d,b) = (a_d b^{-1}, b c_d^{-1})$.
We can recover $d = a_d c_d^{-1}$ from $\iota(d,b)$; since $(a_d,c_d)$ depends
only on $d$, we recover $(a_d,c_d)$ thereby.
From $\iota(d,b)$ and $(a_d,c_d)$ we can tell $b$.
Thus, $\iota$ is an injection.
\end{proof}
In particular, we have
\begin{equation}\label{eq:azr}
d(A,A) \leq d(A, A^{-1}) + d(A^{-1},A) = 2 d(A, A^{-1}) .\end{equation}
If $G$ is abelian, then, by \cite{Ru}, Thm.\ 2,
\begin{equation}\label{eq:ah}
d(A,A^{-1}) \leq 3 d(A,A) .
\end{equation}
This need not hold if $G$ is not abelian: if $A$ is a coset $g H$ of a large 
non-normal subgroup $H\subset G$, we have $|A A^{-1}| = |H| = |A|$, but
$|A A| = |H g H|$ may be much larger than $|A|$, and thus $d(A,A^{-1})$ is
unbounded while $d(A,A) = 0$. 

Another peculiarity of the abelian case is that, if $A \cdot \dotsb \cdot A$
is large, then $A\cdot A$ must be large. If $G$ is not abelian, and
$A$ is of the form $H \cup \{g\}$, where $H$ is a large subgroup
of $G$, then $|A\cdot A| \leq 3 |H| + 1 < 3 |A|$, while 
$A\cdot A\cdot A$ contains $H g H$, and thus may be very large. 
However, the following auxiliary result does hold even for $G$ non-abelian.
\begin{lem}\label{lem:furcht}
Let $n>2$ be an integer. Let $A$ be a finite subset of a group $G$. 
Suppose that
\[|A_{n}| > c |A|^{1 + \epsilon} .\]
for some $c>0$, $\epsilon>0$. Then
\[|A\cdot A \cdot A| > c' |A|^{1 + \epsilon'},\]
where $c'>0$, $\epsilon'>0$ depend only on $c$, $\epsilon$ and $n$.
\end{lem}
\begin{proof}
By (\ref{eq:escolt}),
\[\frac{|A_{n-2} A_2|}{|A|} \leq
\frac{|A_{n-2}\cdot A^{-1}|}{|A|} 
\frac{|A \cdot A_2|}{|A|} \leq \frac{|A_{n-1}|}{|A|} \frac{|A_3|}{|A|}
.\] Proceeding by induction on $n$, we obtain
that
\[\frac{|A_n|}{|A|} \leq \left(\frac{|A_3|}{|A|}\right)^{n-2} .\]
It remains to bound $|A_3|/|A|$ from above by a power of $|A\cdot A\cdot A|/|A|$.
Again by (\ref{eq:escolt}),
\begin{equation}\label{eq:kanad}\begin{aligned}
|A A A^{-1}| |A| &= |A A A^{-1}| |A^{-1}| \leq |A A A| |A^{-1} A^{-1}| \leq
|A A A|^2\\
|A A^{-1} A| |A| &\leq |A A^{-1} A^{-1}| |A A| = |A A A^{-1}| |A A| \leq
|A A A^{-1}| |A A A| .\end{aligned}\end{equation}
Bound $|A A^{-1} A^{-1}|, |A^{-1} A A|,\dotsc, |A^{-1} A^{-1} A^{-1}|$
in terms of $|A A A|$ and $|A|$
by reducing them to either case of (\ref{eq:kanad}): take inverses and
replace $A$ by $A^{-1}$ as needed.
\end{proof}
\subsection{Regularity}
The following is a special case of the Gowers-Balog-Szemer\'edi theorem.
\begin{thm}\label{thm:pope}
Let $A$ be a finite subset of an additive abelian group.
Let $S$ be a subset of $A\times A$ with cardinality
$|S| \geq |A|^2/K$. Suppose we have the bound
\[|\{a+b : (a,b)\in S\}| \leq K |A| .\]
Then there is a subset $A'$ of $A$ such that
$|A'| \geq c K^{-C} |A|$ and
\[|A' + A'| \leq C K^C |A| ,\]
where $c>0$ and $C>0$ are absolute.
\end{thm}
\begin{proof}
By
\cite{Gow}, Prop.\ 12, with $B = A$, there are sets
$A', B'\subset A$ such that $|A'|,|B'| \geq c K^{-C} |A|$
and $|A' - B'|\leq C K^C |A|$. By the pigeonhole principle,
there is a $z$ such that $a - b = z$ for at least
$C^{-1} c^2 K^{-3 C} |A|$ pairs $(a,b) \in A' \times B'$. Thus,
$|V| \geq
C^{-1} c^2 K^{-3 C} |A|$, where we define $V = A'\cap (B'+z)$.
At the same time, $V - V \subset (A' - B') - z$, and so
$|V- V| \leq C K^C |A|$. By (\ref{eq:ah}), $d(V,-V)\leq 3 d(V,V)$, 
and
so $|V + V| \leq \frac{C^6}{c^6} K^{12 C} |V|$. We redefine $A'$ to be $V$
and are done.
\end{proof}
\subsection{Sum-product estimates in finite fields}
\subsubsection{Estimates for small sets}\label{subs:smal}
It is a simple matter to generalize the main result in \cite{Ko} to
finite fields other than $\mathbb{F}_p$.
\begin{thm}\label{thm:bk}
Let $q = p^{\alpha}$ be a prime power. Let $\delta>0$ be given.
 Then, for any $A\subset \mathbb{F}_q^*$
with $C < |A| <p^{1 - \delta}$, we have
\[\max(|A\cdot A|,|A+A|) > |A|^{1 + \epsilon},\]
where $C>0$ and $\epsilon>0$ depend only on $\delta$.
\end{thm}
Explicit values of $C$ and $\epsilon$ can be computed for any 
given $\delta>0$.
\begin{proof}
The proofs of \cite{HBK}, Lem.\ 5, \cite{Ko}, Lem.\ 5, and \cite{Ko}, Thm.\ 2,
work for any finite field $\mathbb{F}_q^*$ without any changes. (In the
statements of \cite{Ko}, Lem.\ 5 and Thm.\ 3, the conditions
$|A| < \sqrt{|F|}$ and $|B| < \sqrt{|F|}$ need to be replaced by
$|A| < \sqrt{p}$ and $|B| < \sqrt{p}$.) For the range $|A| \geq p^{1/2}$,
use \cite{BKT}, Thm.\ 4.3.
\end{proof}
Note the condition $|A| < p^{1-\delta}$ in Thm.\ \ref{thm:bk}, where one
might expect $|A| < q^{1-\delta}$. A subset $A$ of $\mathbb{F}_q^*$ may be
of size about $p$ and fail to grow larger under multiplication by
itself: take, for instance, $A=
(\mathbb{F}_p)^*$, viewed as a subset of $\mathbb{F}_q^*$. One can prove
a version of Thm.\ \ref{thm:bk} in the range $p^{1-\delta} \leq A <
q^{1-\delta}$ (see \cite{BKT}, Thm.\ 4.3), but we will not need to work
in such a range. Hence also the condition $|A| < p^{1-\delta}$ in Prop.
\ref{prop:amtar} and Prop.\ \ref{prop:corz}.
\subsubsection{Estimates for large sets}\label{subs:larg}
\begin{lem}\label{lem:sorge}
Let $p$ be a prime, $A$ a subset of $\mathbb{F}_p$, $S$ a subset
of $\mathbb{F}_p^*$. Then there is an element $\xi \in S$ such that
\[|A + \xi A| \geq \left(\frac{1}{p} + \frac{1}{|S| |A|^2 /p}\right)^{-1}
\geq
\frac{1}{2} \min\left(p, \frac{|S| |A|^2}{p}\right) .\]
Furthermore, for every $c\in (0,1\rbrack$, there are at least
$(1-c) |S|$ elements $\xi \in S$ such that
\[|A + \xi A| \geq c
\left(\frac{1}{p} + \frac{1}{|S| |A|^2 /p}\right)^{-1} .\]
\end{lem}
Cf. \cite{Ko}, Lem.\ 2, which is stronger when $|A| < p^{1/2}$.
\begin{proof}
Let us take Fourier transforms and proceed as in the beginning
of the proof of Thm.\ 6 in \cite{BGK}:
\[\begin{aligned}
p\cdot 
\sum_{\xi\in S} |A\ast \xi A|_2^2 &= \sum_{\xi\in S} 
|\widehat{A\ast \xi A}|_2^2 = \sum_{\xi\in S} |\hat{A} \cdot
\widehat{\xi A}|_2^2 = \sum_{\xi\in S} \sum_{x \in \mathbb{F}_p}
|\hat{A}(x) \widehat{A}(\xi x)|^2 \\
&\leq |S| |\hat{A}(0)|^4 + \sum_{x\in \mathbb{F}_p^*}
\sum_{y\in \mathbb{F}_p^*} |\hat{A}(x) \hat{A}(y)|^2 
= |S| |A|^4 +  \left(\sum_{x\in \mathbb{F}_p^*}
|\hat{A}(x)|^2\right)^2
\\ &= |S| |A|^4 + p^2 (|A|_2^2)^2 =
|S| |A|^4 + p^2 |A|^2 .\end{aligned}\]
Hence, there is an element $\xi_0 \in S$ such that
\[|A\ast \xi_0 A|_2^2 \leq 
\left(\frac{|A|^4}{p} + \frac{p |A|^2}{|S|} \right),\]
and for every $c\in (0,1\rbrack$, 
there are at least $(1 - c) |S|$ elements $\xi \in S$ such that
\[|A\ast \xi A|_2^2 \leq 
\frac{1}{c} \left(\frac{|A|^4}{p} + \frac{p |A|^2}{|S|} \right),\]
By Cauchy's inequality,
\[|A\ast \xi A|_1^2 \leq |A + \xi A| \cdot |A\ast \xi A|_2^2 .\]
As $|A \ast \chi A|_1 = |A|^2$ for every $\chi\in 
\mathbb{F}_p^*$, we obtain that
\[|A + \xi_0 A| \geq \frac{|A \ast \xi_0 A|_1^2}{|A\ast \xi_0 A|_2^2} \geq
\frac{|A|^4}{\frac{|A|^4}{p} + \frac{p |A|^2}{|S|}} =
\left(\frac{1}{p} + \frac{1}{|S| |A|^2 /p}\right)^{-1} \]
for at least one $\xi_0\in S$, and
\[|A + \xi A| \geq \frac{|A \ast \xi A|_1^2}{|A\ast \xi A|_2^2} \geq
\frac{c |A|^4}{\frac{|A|^4}{p} + \frac{p |A|^2}{|S|}} =
c \left(\frac{1}{p} + \frac{1}{|S| |A|^2 /p}\right)^{-1} \]
for at least $(1-c) |S|$ elements $\xi\in S$.
\end{proof}
\section{Expanding functions on $\mathbb{F}_q$}\label{sec:exf}
Let $f$ be a fairly unexceptional
 polynomial on $x$ and $y$ (or on $x$, $x^{-1}$,
$y$ and $y^{-1}$). It is natural to expect a result of the following type
to hold: for every
$\delta>0$ and some $r$, $\epsilon>0$ and $C>0$ depending only on $\delta$,
every set $A\subset \mathbb{F}_p$
 with $C < |A| < p^{1-\delta}$ must fulfill
$|f(A_r,A_r)| > |A|^{1 + \epsilon}$. The work in \cite{BKT} and \cite{Ko}
amounts to such a result for $f(x,y) = x + y$. We will now see how to
derive therefrom a result of the same type for some other choices of
$f(x,y)$.
\begin{prop}\label{prop:amtar}
Let $q = p^{\alpha}$ be a prime power. Let $\delta>0$ be given.
 Then, for any $A\subset \mathbb{F}_q^*$
with $C < |A| <p^{1 - \delta}$, we have
\[|\{(x + x^{-1}) \cdot
(y + y^{-1}) : x,y\in A_2\} | > |A|^{1 + \epsilon},\]
where $C>0$ and $\epsilon>0$ depend only $\delta$.
\end{prop}
\begin{proof}
Let $w(x) = x + x^{-1}$.
Suppose $|\{w(x) w(y) : x,y\in A_2\}| \leq
|A|^{1 + \epsilon}$. It follows directly that $|A_2| \leq \frac{1}{2}
|A|^{1 + \epsilon}$.
Since $w(x) w(y)
= w(x y) + w(x y^{-1})$,
and the cardinality of
$S = \{(w(x y), w(x y^{-1})) : x,y\in A\}$
is at least $|A|^2/16$, we may apply Thm.\ \ref{thm:pope}, and obtain
that there is an $A'\subset A_2$ (which may be taken to be closed under
inversion) such that $|A'| > c' |A|^{1 - C' \epsilon}$ and
$|w(A') + w(A')| < C' |A|^{1 + C' \epsilon}$. At the same time, we have
$|w(A') w(A')| \leq |w(A_2) w(A_2)| \leq |A|^{1 + \epsilon}$.
By Thm.\ \ref{thm:bk}, we have a contradiction, provided that 
$\epsilon$ is small enough and $C$ is large enough.
\end{proof}

\begin{lem}\label{lem:bet}
Let $A$ and $B$ be subsets of a group $G$. Then $A$ can be
covered by at most $|A\cdot B|/|B|$ cosets $a_j B_2$ of $B_2$, where
$a_j\in A$.
\end{lem}
This is the non-commutative version of an argument of Ruzsa's (\cite{Rua}).
\begin{proof}
Let $\{a_1, a_2, \dotsc, a_k\}$ be a maximal subset of $A$ with the
property that the cosets $a_j B$, $1\leq j\leq k$, are all disjoint.
It is clear that $k\leq |A\cdot B|/|B|$. Let $x\in A$. Since $\{a_1,
a_2, \dotsc, a_k\}$ is maximal, there is a $j$ such that $a_j B \cap
x B$ is non-empty. Then $x\in a_j B B^{-1} \subset a_j B_2$. Thus,
the sets $a_j B_2$ cover $A$.
\end{proof}
\begin{prop}\label{prop:corz}
Let $q = p^{\alpha}$ be a prime power.
 Let $\delta>0$ and
$a_1,a_2\in \mathbb{F}_q^*$ be given. Then, for any $A\subset
\mathbb{F}_q^*$ with $C < |A| < p^{1-\delta}$,
\[
|\{a_1 (x y + x^{-1} y^{-1}) +
       a_2 (x^{-1} y + x y^{-1}) : x,y\in A_{20}\}|> |A|^{1 + \epsilon},\]
where $C>0$ and $\epsilon>0$ depend only on $\delta$.
\end{prop}
\begin{proof}
By Lemma \ref{lem:bet}, we may cover $A_4$ with at most $|A_4 \cdot
A^2|/|A^2|$ cosets $a_1 A_2^2, \dotsc, a_k A_2^2$ of $A_2^2$, where
$a_j \in A_4$. Given $x, y\in A_2$ such that $x y\in a_j A_2^2$, we
know that $x y^{-1} = (x y) y^{-2} \in a_j A_4^2$. By Proposition
\ref{prop:amtar} and the pigeonhole principle, there is an index $j$
such that
\begin{equation}\label{eq:cute}
|\{(r + r^{-1}) + (s + s^{-1}) : r,s \in a_j A_4^2\}| >
 \frac{|A|^{1 + \epsilon}}{|A_4 \cdot A^2|/|A^2|} .\end{equation}
Since $|A_4 \cdot A^2|/|A^2| \leq 2|A_6|/|A|$, we have either
$2|A_6|
> |A|^{1 + \epsilon/4}$ or \[\frac{|A|^{1 + \epsilon}}{|A_4 \cdot
A^2|/|A^2|} > |A|^{1 + 3 \epsilon/4}.\] In the former case, we are
already done. So, let us assume $2|A_6| \leq |A|^{1 + \epsilon/4}$.

Write $B = a_j A_4^2\subset A_{12}$. Since
 $|B| \leq |A_4| \leq |A|^{1 + \epsilon/4}$, inequality
(\ref{eq:cute}) implies that \[d_+(w(B),-w(B)) \geq
\frac{\epsilon}{2} \log |A|.\] By (\ref{eq:ah}), we obtain that
\[d_+(w(B),w(B)) \geq \frac{\epsilon}{6} \log |A|.\] Then, by the
triangle inequality (\ref{eq:trian}), \[d_+(a_1 w(B), - a_2 w(B))
\geq \frac{1}{2} d_+(w(B),w(B)) \geq \frac{\epsilon}{12} \log |A|.\]
In other words,
\begin{equation}\label{eq:steme}
|\{ a_1 (r + r^{-1}) + a_2 (s + s^{-1}) : r, s\in B\}| \geq
|B| |A|^{\epsilon/12} \geq \frac{1}{2} |A|^{1 + \epsilon/12} .\end{equation}

For any $r,s\in B$, the ratio $r/s$ is in $A_4^2 A_4^{-2} \subset
A_8^2$. Let $y\in A_8$ be such that $y^2 = r/s$; define $x = r/y \in
A_{20}$. Then $r = x y$ and $s = x/y$. Therefore
\[\{a_1 (r + r^{-1}) + a_2 (s + s^{-1}) : r, s\in B\} \subset
\{a_1 (x y + x y^{-1}) + a_2 (x y^{-1} + x^{-1} y) : x, y\in A_{20}
\} .
\]
By (\ref{eq:steme}), we are done.
\end{proof}
\section{Traces and growth}\label{sec:gerg}
In \S \ref{subs:noeth} we will see how, if 
$A\subset \SL_2(\mathbb{F}_p)$ fails
to grow, it must commute with itself to a fair extent, so to speak.
The arguments in \S \ref{subs:escap} are familiar from the study of growth
in complex groups. The results in \S \ref{subs:sitra} will follow from
those in \S \ref{subs:noeth} by means of simple combinatorial arguments.
We will be able to prove the main part of the key proposition in 
\S \ref{subs:grosma}, using the results in \S \ref{sec:exf}
and \S \ref{subs:noeth}--\ref{subs:sitra}.
\subsection{Growth and commutativity}\label{subs:noeth}
We will first see that, if a subset $A$ of any group $G$ does not grow
rapidly under multiplication by itself, there must be an element $g$
of $A$ with which many elements of $A$ commute. We shall then use the fact
that, in a linear algebraic group, two elements $h_1$, $h_2$ that commute
with a given $g$ with distinct eigenvalues $\lambda_{g,1}, \dotsc,
\lambda_{g,n}$ must also commute with each other. Since non-unipotent
elements are easy to produce in $\SL_2(K)$ (Lem.\ \ref{lem:crud}), we will
conclude that every given
 subset $A$ of $\SL_2(K)$ either grows rapidly or contains
a large simultaneously diagonalizable subset (Cor.\ \ref{cor:kow}).
\begin{prop}\label{prop:tron}
Let $G$ be a group and $A$ a non-empty finite subset thereof.
Let $\Lambda_A$ be the set of conjugacy classes of $G$
with non-zero intersection with $A$. For $g\in G$,
 let $C_G(g)$ be the centralizer of $g$ in $G$.
 Then there is a $g\in A$ such that
\[|C_G(g) \cap (A^{-1} A) | 
\geq \frac{|\Lambda_A| |A|}{|A \cdot A \cdot A^{-1}|} .\]
\end{prop}
\begin{proof}
Let $g,h_1,h_2 \in A$. If $h_1 g h_1^{-1} = h_2 g h_2^{-1}$, then
$h_2^{-1} h_1\in A^{-1} A$ commutes with $g$. Hence, for any $g\in G$, 
\[|\{h g h^{-1} : h\in A\}| \geq \frac{|A|}{|C_G(g) \cap A^{-1} A|} .\]
Let $\Upsilon \subset A$ be a set of representatives of $\Lambda_A$.
Then
\[|A A A^{-1}| \geq |\{h g h^{-1} : h\in A, g\in \Upsilon\}| \geq
\sum_{g\in \Upsilon} \frac{|A|}{|C_G(g) \cap A^{-1} A|} .\]
If $|C_G(g) \cap (A^{-1} A)|  < 
\frac{|\Lambda_A| |A|}{|A \cdot A \cdot A^{-1}|}$ for every $g\in \Upsilon$,
then 
\[\sum_{g\in \Upsilon} \frac{|A|}{|C_G(g) \cap A^{-1} A|} 
> |\Upsilon| \frac{|A\cdot A\cdot A^{-1}|}{|\Lambda_A|} = 
|A\cdot A\cdot A^{-1}|,\]
and we reach a contradiction.
\end{proof}
\begin{lem}\label{lem:crud} 
Let $K$ be a field. Let $A$ be a finite subset of $\SL_2(K)$ not
contained in any proper subgroup of $\SL_2(K)$. Then
$A_2$ has at least $\frac{1}{4} |A| - 1$ elements with trace
other than $\pm 2$.
\end{lem}
\begin{proof}
Let $g\in A$ be an element of trace $2$
or $-2$ other than $\pm I$. Let $B\subset A$ 
be the set of all elements of $A$
 with trace $\pm 2$ and an eigenvector in common
with $g$. Suppose
$|B| \leq \frac{1}{4} |A| + 3$. Let $h\in A\setminus B$. 
If $h$ has trace $\pm 2$, then either $g h$ 
or $g^{-1} h$ does not. Therefore $A \cup A\cdot A \cup A^{-1} A$ 
has at least
$\frac{1}{3} |A\setminus B| \geq \frac{1}{4} |A| - 1$ elements with trace 
other than
$2$. Suppose now $|B|> \frac{1}{4} |A| + 3$. Let $h$ be an element of $A$
that does not have an eigenvector in common with $g$. Then there are at most
two elements $g'$ of $B$ such that $g' h$ has trace $2$. Hence
$A \cdot A$ has more than $\frac{1}{4} |A| + 1$ 
elements with trace other than $\pm 2$. 
\end{proof}
\begin{cor}\label{cor:kow} 
Let $K$ be a field. Let $A$ be a non-empty finite subset of $\SL_2(K)$ not
contained in any proper subgroup of $\SL_2(K)$. 
Assume $|\Tr(A)|\geq 2$,
$|A|\geq 4$. Then 
there are at least 
$\frac{(|\Tr(A)| - 2) (\frac{1}{4} |A| - 1)}{|A_6|}$
simultaneously diagonalizable matrices in $A_4$.
\end{cor}
\begin{proof}
Let $B$ be the set of elements of 
$A_2$ with trace other than
$\pm 2$.
By Lemma \ref{lem:crud}, $|B| \geq \frac{1}{3} |A| - 1$.
We may apply Prop.\ \ref{prop:tron}, and obtain
that there is a $g\in B$ such that
\[|C_G(g) \cap (B^{-1} B)| \geq 
\frac{|\Lambda_B| |B|}{|B\cdot B \cdot B^{-1}|} \geq
\frac{|\Tr(B)| |B|}{|B\cdot B\cdot B^{-1}|} \geq
\frac{(|\Tr(A)| - 2) (\frac{1}{4} |A| - 1)}{|A_6|}.\]
All elements of $V = C_G(g) \cap (B^{-1} B)$ commute with $g$;
since $\Tr(g)\ne \pm 2$, it follows that,
when $g$ is diagonalized, so is all of $V$.
\end{proof}
\subsection{Escaping from subvarieties}\label{subs:escap}
 The following lemma\footnote{Thanks are
  due to N.\ Anantharaman for pointing out an inaccuracy in a previous version
of this paper, and to both N.\ Anantharaman and E.\ Breuillard for help with the current 
phrasing.} is based closely on \cite[Prop.\ 3.2]{EMO}.
\begin{lem}\label{lem:carbo}
Let $G$ be a group. Consider a linear representation of $G$
on a vector space $V$ over a field $K$. Let $W$ be a union
$W_1 \cup W_2 \cup \dotsc \cup W_n$ of proper subspaces of $V$.

Let $A$ be a subset of $G$; let $\mathscr{O}$ be an $\langle A\rangle$-orbit
in $V$ not contained in $W$. Then there are constants $\eta>0$ and $m$
depending only on $n$ and $\dim V$ such that, for every $x\in \mathscr{O}$,
there are at least $\max(1, \eta |A|)$ elements $g\in A_m$ such that
$g x\notin W$.
\end{lem}
This may be phrased as follows: one can escape from $W$ by the action
of the elements of $A$. One can give stronger and more general statements
of this kind; the spaces $W_n$ could very well be taken to be varieties instead.
However, what we have just stated will do.

\begin{proof}
Let us begin by showing that there are elements $g_1,\dotsc,g_l \in A_{r}$
such that, for every $x\in \mathscr{O}$, at least one of the
$g_i \cdot x$'s is not in $W$. (Here $l$ and $r$ are bounded
in terms of $n$ and $d = \dim V$ alone.) We will proceed by induction
on $(d_W,s_W)$, where $d_W$ is the maximal dimension of the spaces
$W_1, \dotsc, W_n$ (i.e., $d_W = \max_{1\leq j\leq n} \dim(W_j)$)
and $s_W$ is the number of spaces of dimension $d_W$ among
$W_1,\dotsc,W_n$. We shall always pass from $W$ to a union of the form
$W' = W_1' \cup \dotsb \cup W_n'$, where either (a) $d_{W'} < d_W$ or
(b) $d_{W'} = d_W$ and $s_{W'} < s_W$. The base case of the inductive
process will be $(d_W,s_W) = (0,0)$.


Let $W_+$ be the union of subspaces $W_j$, $1\leq j\leq n$, of dimension
$d_W$ (the maximal dimension). If $W_+$ and $\mathscr{O}$ are disjoint, we set
$W' = W\setminus W_+$. Suppose otherwise.
Since $\mathscr{O}$ is not contained in $W_+$, we can find
$x_0\in W_+\cap \mathscr{O}$, $g\in A \cup A^{-1}$ such that
$g x_0 \notin W_+$. 
Hence the set of subspaces of maximal dimension in $W$ is not the same
as the set of subspaces of maximal dimension in $W'$.
It follows that $W' = g W \cap W$ 
does not contain $W_+$, and thus has fewer subspaces $W_j'$ of dimension
$d_W$ (the maximal dimension) than $W$ has. 

We have thus passed from $W$ to $W'$, where either (a) $d_W' < d_W$ or
(b) $d_W' = d_W$ and $s_W' < s_W$. By the inductive hypothesis, we 
already know that there are
$g_1',\dotsc,g_{l'}' \in A_{r'}$
such that, for every $x\in \mathscr{O}$, at least one of the
$g_i' \cdot x$'s is not in $W'$. (Here $l'$ and $r'$ are bounded
in terms of $n'$ and $d = \dim V$ alone; the number $n'$ of subspaces
$W_1', W_2',\dotsc , W_{n'}'$ is bounded by $n^2$.) 
Since at least one of the $g_i' \cdot x$'s
is not in $W' = g W\cap W$, either one of the $g_i'\cdot x$'s is not
in $W$ or one of the $g_i'\cdot x$'s is not in $g W$, i.e.,
one of the $g^{-1} g_i' \cdot x$'s is not in $W$. Set
\[\begin{aligned}
g_1 &= g_1',\; g_2 = g_2',\; \dotsc,\; g_l = g_l'\\
g_{l+1} &= g^{-1} g_1',\; g_{l+2} = g^{-1} g_2',\; \dotsc,\; g_{2 l} = g^{-1}
g_l',\;\;\;\;\; l' = 2 l .\end{aligned}\]
(As can be seen, $g_i \in A_r$, where $r = r'+1$.)
We conclude that, for every $x\in \mathscr{O}$, at least one of
the $g_i \cdot x$'s is not in $W$.

The rest is easy: for each $x\in \mathscr{O}$ and each $g\in A$,
at least one of the elements $g_i g \cdot x$, $1\leq i\leq l$ 
($g_i\in A_r$)
will not be in $W$.
Each possible $g_i g$ can occur for at most $l$ different elements $g\in A$;
thus, there are at least $\min(1,|A|/l)$ elements $h = g_i g$ 
of $A_{r+1}$ such that $h x \notin W$.
\end{proof}
We derive some immediate consequences.

\begin{cor}\label{cor:rats}
Let $K$ be a field.
Let $A$ be a finite subset of $\SL_2(K)$ not
contained in any proper subgroup of $\SL_2(K)$.
If $|K|>3$, the following holds\/{\rm :}\/
for any basis
$\{v_1, v_2\}$ of $\overline{K}^2${\rm ,} there is a $g\in A_{k}$ such that
 $g v_i \ne \lambda v_j$ for all choices of $\lambda\in \overline{K}${\rm ,}
$i,j\in \{1,2\}${\rm ,} where $k$ is an absolute constant.
\end{cor}
\begin{proof}
Consider $G = \SL_2(K)$ and its natural action on the vector
space $V = M_2(\overline{K})$ of $2$-by-$2$ matrices.
Let $W$ be the subset of $V$ consisting of all $h\in V$ such that
$h v_i = v_j$ for some $i,j\in \{1,2\}$.
Let $x$ be the identity in $M_2(\overline{K})$.
Apply Lemma \ref{lem:carbo}.

Before Lemma \ref{lem:carbo} can be applied, we must verify\footnote{
Thanks to O. Dinai for the counting argument about to be used.}
 that the orbit $\mathscr{O} = \SL_2(K)$ of $x$ is not contained in $W$.
Let $G_{i,j}$ be the set of matrices $g$ in $\SL_2(K)$ such that
$g v_i$ is a multiple of $v_j$. Since $W(K) \cap \mathscr{O} 
= G_{1,1} \cup G_{1,2} \cup G_{2,1}
\cup G_{2,2}$, we would like to bound $|G_{i,j}|$. Let $g\in G_{i,j}$.
Choose a vector $v\in
K^2$ (say $v = (1,0)$ or $v=(0,1)$) that is not a multiple of $v_i$. 
It is clear that $g v$ and $g v_i$ determine $g$. At the same time, 
we already know that $g v_i = \lambda v_j$, and, if $g v$ is fixed,
 two different values of $\lambda$ determine two matrices $g$ with
different determinants; in particular, at most one $\lambda\in
\overline{K}$ gives us
a $g\in \SL_2(K)$. Thus $g v$ actually determines $g$. Since $g v$ must
be non-zero and lie in $K^2$, we conclude that $|G_{i,j}| \leq |K|^2 - 1$.

The sets $G_{1,1}$ and $G_{2,2}$ intersect at the identity. Thus,
$|W(K) \cap \mathscr{O}| \leq 4 (|K|^2 - 1) - 1$. Since $|\SL_2(K)| = |K| \cdot (|K|^2 - 1)$,
it is enough to assume $|K|\geq 4$ to conclude that $|W(K)\cap \mathscr{O}| 
< |\SL_2(K)|$.
In particular, for $|K|\geq 4$, the set $\mathscr{O} = \SL_2(K)$ is not
contained in $W$. We are entitled to apply Lemma \ref{lem:carbo}, after all.
\end{proof}

\begin{cor}\label{cor:kot}
Let $K$ be a field. Let $A$ be a finite subset of $\SL_2(K)$ not
contained in any proper subgroup of $\SL_2(K)$.
Then there
are absolute constants $k,c>0$ such that{\rm ,} 
given any two non-zero
vectors
$v_1, v_2\in \overline{K}^2${\rm ,}
\[|A_k\setminus (H_{v_1}\cup H_{v_2})|
> c |A| ,\]
where $H_v = \{g\in \SL_2(K): \text{$v$ is an eigenvector of $g$}\}$.
\end{cor}
\begin{proof}
Consider $G = \SL_2(K)$ and its natural action on $V = M_2(\overline{K})$.
Let $W = H_{v_1}' \cup H_{v_2}'$, where
$H_v' = \{g\in M_2(\overline{K}): \text{$v$ is an eigenvector of $g$}\}$.
Let $x = I$. 

Before we apply Lemma \ref{lem:carbo}, we need to check that
$\SL_2(K)$ is not contained in $W(K)$. Since the matrices 
$\left(\begin{matrix}1 & 1\\0 & 1\end{matrix}\right)$,
$\left(\begin{matrix}1 & 0\\1 & 1\end{matrix}\right)$ and
$\left(\begin{matrix}0 & 1\\-1 & 0\end{matrix}\right)$ share no
eigenvectors, there is no pair of eigenvectors $v_1$, $v_2$ such that
each of three matrices has at least one of $v_1$, $v_2$ as an eigenvector.
Thus $\SL_2(K) \not\subset W(K)$. Now apply Lemma~\ref{lem:carbo}.
\end{proof}

Lemma \ref{lem:crud} could be derived from Lemma \ref{lem:carbo} as well,
but, since the proof of Lemma \ref{lem:crud} is simple as it is,
 we will not bother.
\subsection{Size from trace size}\label{subs:sitra}
Given a large set $V$ of diagonal matrices and a matrix $g\notin V$
with only non-zero entries, one can multiply $V$ and $g$ to obtain
at least $\gg |V|^3$ different matrices.
\begin{lem}\label{lem:chich}
Let $K$ be a field.
Let $V\subset \SL_2(K)$ be a
 finite set of simultaneously diagonalizable matrices; call their
common eigenvectors $v_1$ and $v_2$. Let
$g\in \SL_2(K)$ be such that $g v_i \ne \lambda v_j$ for any
$\lambda\in \overline{K}$, $i,j\in \{1,2\}$. Then 
\[|V g V g^{-1} V| \geq \frac{1}{2} \left(\frac{1}{4} |V| - 5\right) |V|^2 .
\]
\end{lem}
\begin{proof}
Diagonalize $V$, conjugating by an element of $\SL_2(\overline{K})$ if
necessary. Write $g = \left(\begin{array}{cc} a & b\\
c & d\end{array}\right)$. By assumption, $a b c d\ne 0$.
Then
\begin{equation}\label{eq:mandru}
 g \left(\begin{array}{cc} r & 0\\ 0 & r^{-1}\end{array}\right) g^{-1}
= \left(\begin{array}{cc} r a d - r^{-1} b c &
(r^{-1} - r) a b \\ (r - r^{-1}) c d & r^{-1} a d - r b c\end{array}\right) ,
\end{equation}
the product of whose upper-right and lower-left entries is
$-(r-r^{-1})^2 a b c d$. 
The map $r\mapsto -(r-r^{-1})^2 a b c d$ cannot send more than
$4$ distinct elements of $K^*$ 
to the same element of $K$. Thus, the set 
$\{ h_{12} h_{21} : h\in g V g^{-1} \}$ has cardinality at least $|V|/4$.
The upper-left and lower-right entries of the matrix in the right-hand side
of (\ref{eq:mandru}) can be both equal to $0$ only if $r^2 - r^{-2} = 0$,
and that can happen for at most $4$ values of $r$. Let 
$U = \{h\in g V g^{-1} : (h_{11} h_{12} h_{21} \ne 0) 
\wedge (h_{22} h_{12} h_{21}\ne 0)\}$; we have that 
$|\{h_{12} h_{21} : h\in U\}| \geq \frac{1}{4} |V| - 5$.

Let $h\in U$ be fixed. Define
\[f_h(s,t) = \left(\begin{array}{cc} s & 0\\ 0 & s^{-1}\end{array}\right) 
\left(\begin{array}{cc} h_{11} & h_{12}\\ h_{21} & h_{22}\end{array}\right) 
\left(\begin{array}{cc} t & 0\\ 0 & t^{-1}\end{array}\right)
= \left(\begin{array}{cc} s t h_{11} & s t^{-1} h_{12}\\ 
s^{-1} t h_{21} & s^{-1} t^{-1} h_{22}\end{array}\right).\]
The product of the upper-right and lower-left entries of $f_h(s,t)$ is
$h_{12} h_{21}$, which is independent of $s$ and $t$. Since
$h\in U$, we may recover
$s^2$, $t^2$ and $s t$ from $h$ and $f_h(s,t)$. Thus, for $h$ fixed,
there cannot be more than two pairs $(s,t)$ sharing the same value of
$f_h(s,t)$. For each element of
$\{ h_{12} h_{21} : h\in U \}$, choose an $h$ corresponding to it;
let $s$ and $t$ vary. We obtain at least $\frac{1}{2} 
|\{h_{12} h_{21} : h\in U\}| |V|^2$
different values of $f_h(s,t) \in V g V g^{-1} V$.
 We conclude that
$\{ V g V g^{-1} V\}$ has cardinality at least 
$\frac{1}{2} |\{h_{12} h_{21} : h\in U\}| |V|^2 = 
\frac{1}{2} (\frac{1}{4} |V| - 5) |V|^2$.
\end{proof}
We will now use Cor.\ \ref{cor:kow}, Cor.\ \ref{cor:rats}
and Lem.\ \ref{lem:chich} to show that, unless $A$ grows
substantially under multiplication by itself, the cardinality
of $A_k$ cannot be much smaller than the cube of the
cardinality of the set of traces $\Tr(A)$ of $A$. 
\begin{prop}\label{prop:andu}
Let $K$ be a field. Let $A$ be a finite subset of $\SL_2(K)$ not
contained in any proper subgroup of $\SL_2(K)$. Assume $|\Tr(A)|\geq 2$,
$|A|\geq 4$ and $|K|>3$.
Then 
\[|A_k| \geq 
\frac{1}{2} \left(\frac{1}{4}
\frac{(|\Tr(A)|-2) (\frac{1}{4} |A| - 1)}{|A_6|} - 5\right)
\left(\frac{(|\Tr(A)|-2) (\frac{1}{4} |A| - 1)}{|A_6|}\right)^2
,\]
where $k$ is an absolute constant. 
\end{prop}
\begin{proof}
By Cor.\ \ref{cor:kow}, there is a simultaneously diagonalizable subset 
$V \subset A_4$ with 
$|V| \geq \frac{(|\Tr(A)|-2) (\frac{1}{4} |A| - 1)}{|A_6|}$; call
its common eigenvectors $v_1$ and $v_2$.
Since $A$ is not contained in any proper subgroup of
$\SL_2(K)$, Cor.\ \ref{cor:rats} yields a $g\in A_k$ such that
$g v_i \ne \lambda v_j$ for all $\lambda\in K$, $i,j\in \{1,2\}$.
Hence, by Lemma \ref{lem:chich},
$|V g V g^{-1} V| \geq 
\frac{1}{2} \left(\frac{1}{4} |V| - 5\right) |V|^2$.
\end{proof}
We must now prove that, unless $A$ grows substantially when
multiplied by itself, the cardinality of $\Tr(A_k)$ cannot
be much smaller than the cube root of the cardinality of $A$. A
preparatory lemma is needed. Like Lem.\ \ref{lem:chich}, it
is of a very simple type -- the cardinality of a set is bounded
from below by virtue of its being contained the image of a map 
that has a large enough domain and is not too far from being injective.
\begin{lem}\label{lem:funn}
Let $K$ be a field. Let $A$ be a finite subset of $\SL_2(K)$. Write
the matrices in $\SL_2(K)$ with respect to a basis $\{v_1, v_2\}$
of $\overline{K}^2$.
Suppose $g_{12} g_{21} \ne 0$ for every
$g\in A$. Then
\[|\Tr(A A^{-1})| \geq 
\frac{|A|}{2\cdot |\{(g_{11}, g_{22}) : g\in A\}|} .\]
\end{lem}
\begin{proof}
Let $D = \{(g_{11}, g_{22}) : g\in A\}$. Consider any two
distinct
$g, g'\in B$ with $g_{11} = g_{11}'$, $g_{22} = g_{22}'$.
Then $g g'^{-1}$ has trace
\[\Tr(g g'^{-1}) = g_{11} g_{22}' + g_{22} g_{11}'
- g_{12} g_{21}' - g_{21} \left(\frac{g_{11}' g_{22}' - 1}{g_{21}'}
\right) .\]
Thus, given $g\in B$, 
there can be at most two $g'\in B$ with
$g_{11} = g_{11}'$, $g_{22} = g_{22}'$ such that $\Tr(g g'^{-1})$
is equal to a given value. Choose $g$ such that
$|\{g'\in B: g'_{11} = g_{11}, g'_{22} = g_{22}\}|$ is maximal.
\end{proof}
\begin{prop}\label{prop:unda}
Let $K$ be a field. Let $A$ be a finite subset of $\SL_2(K)$
not contained in any proper subgroup of $\SL_2(K)$. Then
\[|\Tr(A_k)| \geq c |A|^{1/3},\]
where $k$ and $c>0$ are absolute constants.
\end{prop}
\begin{proof}
If $A$ has an element of trace other than $\pm 2$, let $h$ be one such
element. Otherwise, choose any $g_1\in A$ other than
$\pm I$, and any $g_2\in A$ not in the
unique Borel subgroup in which $g_1$, being parabolic, lies; then either
$g_1 g_2\in A\cdot A$ or $g_1^{-1} g_2 \in A^{-1} A$ has trace
$\ne \pm 2$; choose $h\in A_2$, $\tr(h) \ne \pm 2$, to be one of the two.
From now on, write all matrices
with respect to the two eigenvectors $v_1$, $v_2$ of $h$.
We denote by $r$ and $r^{-1}$ the two eigenvalues of $h$.

By Cor.\ \ref{cor:kot}, $|X| \geq c |A|$, where
$X = A_{k_0}\setminus (H_{v_1} \cup H_{v_2})$ and
$k$, $c>0$ are absolute constants. 
Lemma \ref{lem:funn} now implies that
\begin{equation}\label{eq:unghu}
|\Tr(A_{2 k_0})| \geq |\Tr(X X^{-1})| \geq 
\frac{|X|}{2\cdot |\{(g_{11}, g_{22}) : g\in X\}|} .\end{equation}

For $t\in K$, let $D_t =
|\{(g_{11}, g_{22}): g_{11} + g_{22} = t, g\in X\}|$. Let
$t\in K$ be such that $|D_t|$ is maximal. For any 
$(a,d) \in D_t$, we have
$r a + r^{-1} d = (r - r^{-1}) a + r^{-1} t$. Thus, for any two
distinct pairs $(a,d), (a',d')\in D_t$, the two values
$r a + r^{-1} d$, $r a' + r^{-1} d'$ must be distinct. 
Thus 
\[|\Tr(A_{k_0+2})| \geq |\Tr(h X)| \geq |D_t| \geq \frac{
|\{(g_{11}, g_{22}) : g\in X\}|}{|\Tr(X)|} .\]
Multiplying by (\ref{eq:unghu}), we obtain
\[|\Tr(A_{k_0+2})| |\Tr(A_{2 k_0})|  \geq \frac{|X|}{2 |\Tr(X)|} ,\]
and so $|\Tr(A_{2 k_0})|^3 \geq 
|\Tr(A_{k_0+2})| |\Tr(A_{2 k_0})| |\Tr(X)| \geq \frac{1}{2} |X|$,
where we assume, as we may, that $k_0\geq 2$.
Hence
 \[|\Tr(A_{2 k_0})|\geq \left(\frac{1}{2} |X|\right)^{1/3} \geq 
\frac{c_0^{1/3}}{2^{1/3}} |A|^{1/3} .\]
\end{proof} 
\subsection{Growth of small sets}\label{subs:grosma}
The statements in the section up to now reduce the main problem
to a question in $\mathbb{F}_{p^2}$, and that question can be answered
using the results in \S \ref{sec:exf}.
\begin{proof}[Proof of part (\ref{it:wark}) of the key proposition]
We may assume that $p$ is larger than an absolute constant;
otherwise we may make (\ref{eq:solt}) true simply by adjusting
the constant $c$ therein. By the same token, we may assume that
$|A|$ is larger than an absolute constant.

By Proposition \ref{prop:unda}, 
$|\Tr(A_{k_0})| \geq c_0 |A|^{1/3}$, where $k_0$ and $c_0$
are absolute constants. 
As we said, we may assume that
$|A| \geq \max((4/c_0)^3,8)$.
Thus, by
Cor.\ \ref{cor:kow}, there are at least
\[\frac{(c_0 |A|^{1/3} - 2) (\frac{1}{4} |A_{k_0}| - 1)}{|A_{6 k_0}|}
\geq \frac{c_0 |A|^{1/3} |A_{k_0}|}{16 |A_{6 k_0}|}\]
simultaneously diagonalizable matrices in $A_{4 k_0}$; 
denote by $V$ the set
of the eigenvalues of $\lceil \frac{c_0 |A|^{1/3} |A_{k_0}|}{16 |A_{6 k_0}|}
\rceil$
such matrices. 
Since we may 
assume that $c_0<1$, we have $|V| < |A|^{1/3} < p^{1 - \delta/3}$.
We also take for granted that $|A_{6 k_0}| < |A|^{7/6}$; otherwise, by Lem.\
\ref{lem:furcht}, we are already done. 
Thus $|V| > \frac{c_0}{16} |A|^{1/6}$, and so, given
a $C$ depending only on $\delta$, we may assume that $|V|>C$
by adjusting the constant $c$ in (\ref{eq:solt}) accordingly.

By Corollary \ref{cor:rats}, there is a matrix
$\left(\begin{array}{cc} a & b\\ c & d\end{array}\right) \in A_{k_1}$ such that
$a b c d \ne 0$, where $k_1$ is an absolute constant.
 Now, for any scalars $x,y$, the trace of
\[\left(\begin{array}{cc} x & 0\\ 0 & x^{-1}\end{array}\right)
\left(\begin{array}{cc} a & b\\ c & d\end{array}\right)
\left(\begin{array}{cc} y & 0\\ 0 & y^{-1}\end{array}\right)
\left(\begin{array}{cc} d & -b\\ -c & a\end{array}\right)\] 
is $a d (x y + x^{-1} y^{-1}) - b c (x^{-1} y + x y^{-1})$. Letting
$x$, $y$ range on all of $V$, we see
that $\tr(A_{160 k_0 + 2 k_1}) = 
\tr(A_{20\cdot 4 k_0 + k_1 + 20\cdot 4 k_0 + k_1}) \supset 
\{a d (x y + x^{-1} y^{-1}) - 
b c (x^{-1} y + x y^{-1}) : x,y\in V_{20}\}$. Now we apply Prop.\ 
\ref{prop:corz} with $q=p^2$, and obtain that
\[|\tr(A_{160 k_0 + 2 k_1})| > |V|^{1 + \epsilon},\]
where $\epsilon>0$ depends only on $\delta$.
Here we have assumed, as we may, 
that $|V|>C$, where $C$ is the constant
in the statement of Prop.\ \ref{prop:corz}, with $\delta$ equal to one-third
of our $\delta$.

By the same argument as when we took $|V| > \frac{c_0}{16} |A|^{1/6}$,
 we may assume that
\[\frac{|\Tr(A_{160 k_0 + 2 k_1})| |A_{160 k_0 + 2 k_1}|}{|A_{6
(160 k_0 + 2 k_1)}|}\geq 40 .\] 
(Otherwise we are already done.)
We proceed by applying Prop.\ \ref{prop:andu}, and obtain
\[\begin{aligned}
|A_{k_2 (160 k_0 + 2 k_1)}| &\geq \frac{1}{2^{16}}
\frac{|\Tr(A_{160 k_0 + 2 k_1})|^3 |A_{160 k_0 + 2 k_1}|^3}{|A_{6
(160 k_0 + 2 k_1)}|^3} 
> \frac{1}{2^{16}}
\frac{|A_{160 k_0 + 2 k_1}|^3}{|A_{6
(160 k_0 + 2 k_1)}|^3} 
|V|^{3 (1 + \epsilon)} \\ &\geq \frac{1}{2^{16}}
\frac{|A_{160 k_0 + 2 k_1}|^3}{|A_{6
(160 k_0 + 2 k_1)}|^3} 
\frac{c_0^3 |A_{k_0}|^3}{2^{12} |A_{6 k_0}|^3} |A|^{1+\epsilon} \geq
\frac{c_0^3}{2^{28}} \frac{|A|^6}{|A_{6
(160 k_0 + 2 k_1)}|^6} |A|^{1 + \epsilon},\end{aligned}\]
where $k_2$ is an absolute constant. 
Hence, either $|A_{6 (160 k_0 + 2 k_1)}|$ or 
$|A_{k_2 (160 k_0 + 2 k_1)}|$ must be greater than
$\frac{c_0^{3/7}}{16} |A|^{1 + \epsilon/7}$.
By Lemma \ref{lem:furcht}, we are done.
\end{proof}
\section{Generating the whole group}\label{sec:whg}
Since we have proved part (\ref{it:wark}) of the key proposition, we
know how to attain a set of cardinality $p^{3 - \delta}$, $\delta>0$,
by multiplying a given set of generators 
$A$ by itself $(\log(p/|A|))^c$ times. It
remains to show how to produce the 
group $\SL_2(\mathbb{Z}/p \mathbb{Z})$ in a bounded number of steps
from a set almost as large as $\SL_2(\mathbb{Z}/p \mathbb{Z})$ itself.
As might be expected,
instead of the sum-product estimates for small sets (\S \ref{subs:smal}), 
we will use the estimates for large sets (\S \ref{subs:larg}).
We first focus on what happens in the Borel subgroups. 
\begin{lem}\label{lem:attac}
Let $p$ be a prime. Let $H$ be a Borel subgroup of 
$\SL_2(\mathbb{Z}/p \mathbb{Z})$. Let $A\subset H$ be given with
$|A| > 2 p^{5/3} + 1$. Then $A_8$ contains all elements of $H$ with trace $2$.
\end{lem}
\begin{proof}
We may as well assume that $H$ is the set of upper-triangular 
matrices. Define $P_r(A) = 
\left\{x\in \mathbb{Z}/p \mathbb{Z} :
\left(\begin{array}{cc} r & x\\ 0 & r^{-1}\end{array}\right)\in A
\right\}$.
By the pigeonhole principle, there is an $r\in
(\mathbb{Z}/p \mathbb{Z})^*$ such that 
$|P_r(A)| > 2 p^{2/3}$. Let
$\left(\begin{array}{cc} t & u\\ 0 & t^{-1}\end{array}\right)$
be any element of $A$ with $t\ne r$. Then
\[
\left(\begin{array}{cc} t & u\\ 0 & t^{-1}\end{array}\right)
\left(\begin{array}{cc} r & x\\ 0 & r^{-1}\end{array}\right)
\left(\begin{array}{cc} t^{-1} & -u\\ 0 & t\end{array}\right)
\left(\begin{array}{cc} r^{-1} & -x'\\ 0 & r\end{array}\right)\]
equals
\[
\left(\begin{array}{cc} r & t^2 x + (r^{-1} - r) u t\\ 
0 & r^{-1}\end{array}\right)
\left(\begin{array}{cc} r^{-1} & -x'\\ 0 & r\end{array}\right) \\
=
\left(\begin{array}{cc} 1 & r (- x' + t^2 x) + (1 - r^2) u t\\
0 & 1\end{array}\right) .\]
Therefore, $P_1(A A A^{-1} A^{-1})$
is a superset of $r (-P_r(A) + t^2 P_r(A)) + (1 - r^2) u t$. Define
$S = \{t\in (\mathbb{Z}/p \mathbb{Z})^*,\, t\ne r: \exists\, 
u\in \mathbb{Z}/p \mathbb{Z}\; \text{s.t.} 
\left(\begin{array}{cc} t & u\\ 0 & t^{-1}\end{array}
\right) \in A,\, u\in \mathbb{Z}/p \mathbb{Z},\,
 t\ne r\}$. Clearly $|S| > \frac{1}{p} (2 p^{5/3} - p) > p^{2/3}$. By Lemma
\ref{lem:sorge}, there is a $t \in S$ such that
\[
|r (-P_r(A) + t^2 P_r(A)) + (1 - r^2) u t| = |P_r(A) - t^2 P_r(A)|
\geq
\frac{1}{\frac{1}{p} + \frac{p}{\frac{1}{2} |S| |P_r(A)|^2}} >
\frac{1}{\frac{1}{p} + \frac{1}{2 p}} = \frac{2}{3} p .\]
Thus,
\[(r (P_r(A) + t^2 P_r(A)) + (1 - r^2) u t) + 
(r (P_r(A) + t^2 P_r(A)) + (1 - r^2) u t) = \mathbb{Z}/p \mathbb{Z} .\]
It follows that $A A A^{-1} A^{-1} A A A^{-1} A^{-1}$ contains
all matrices $\left(\begin{array}{cc} 1 & x\\ 0 & 1\end{array}
\right)$, $x\in \mathbb{Z}/ p \mathbb{Z}$.
\end{proof}
\begin{proof}[Proof of part (\ref{it:agust}) of the key proposition]
By part (\ref{it:wark}) of the main theorem, we may assume that
$|A| > 6 p^{8/3} > (2 p^{5/3} + 1) ( p + 1)$.
By the pigeonhole principle, there are at least
$(2 p^{5/3} + 1)$ matrices in $A$ with the same lower row
up to multiplication by a scalar in $(\mathbb{Z}/ p \mathbb{Z})^*$; 
the same holds, of course, for the upper row. Thus,
there are at least $2 p^{5/3}+1$ upper-diagonal matrices
and at least $2 p^{5/3}+1$ lower-diagonal matrices in
$C = A A^{-1}$.
By Lemma \ref{lem:attac},
$C_8$ contains all matrices of the form
$\left(\begin{array}{cc} 1 & x\\ 0 &1\end{array}\right)$,
$\left(\begin{array}{cc} 1 & 0\\ y &1\end{array}\right)$,
$x,y\in \mathbb{Z}/ p \mathbb{Z}$. Every element of $\SL_2(\mathbb{Z}/p \mathbb{Z})$
can be written in the form
\[\left(\begin{array}{cc} 1 & 0\\ y &1\end{array}\right)
  \left(\begin{array}{cc} 1 & x\\ 0 &1\end{array}\right)
\left(\begin{array}{cc} 1 & 0\\ y' &1\end{array}\right)
  \left(\begin{array}{cc} 1 & x'\\ 0 &1\end{array}\right),\]
where $x,y,x',y'\in \mathbb{Z}/ p \mathbb{Z}$. Hence
$\SL_2(\mathbb{Z}/p \mathbb{Z}) = C_8 C_8 C_8 C_8 \subset A_{64}$.
\end{proof}
{\em Note added in proof.} A far more elegant proof of part (b)
given part (a) may be obtained by an approach due to Gowers \cite{Gow2}; 
see \cite{NP}. In brief: in the present context,
 it is cleaner and simpler to do Fourier analysis
on $\SL_2(\mathbb{Z}/p\mathbb{Z})$ itself, rather than to prove and use
results based on Fourier analysis over $\mathbb{Z}/p\mathbb{Z}$
(\S \ref{subs:larg}, \S \ref{sec:whg}).
\section{The main theorem and further consequences}\label{sec:conc}
\begin{proof}[Proof of Main Theorem] The statement of the theorem follows
immediately from the key proposition, 
parts (\ref{it:wark}) and (\ref{it:agust}),
when $|A|$ is larger than an absolute constant. 
Since $|A \cup A \cdot A|\geq
|A| + 1$ for any $A$ not a subgroup of 
$\SL_2(\mathbb{Z}/p \mathbb{Z})$, we may increase
the cardinality of $A$ by an absolute constant $C$
 simply by multiplying $A$ by itself $C$ times.
\end{proof}

Let $G$ be a finite group and $A\subset G$ a set of generators of $G$.
Let $\psi$ be a probability distribution on $G$ whose support contains $A$.
We will assume throughout that $\psi$ is symmetric, i.e., 
$\psi(g) = \psi(g^{-1})$ for every $g\in G$.
We define
the transition matrix $T_{\psi}(G,A) = \{\psi(y^{-1} x)\}_{x,y\in G}$.
The largest eigenvalue of $T_{\psi}(G,A)$ is clearly $1$. 

Consider a 
family $\{G_j,A_j\}_{j\in J}$ 
of finite groups $G_j$ and sets of generators $A_j$ of $G_j$
such that $d = |A_j \cup A_j^{-1}|$ is constant. Let
$\psi_j(g) = \frac{1}{d}$ if $g\in A_j\cup A_j^{-1}$ and $\psi_j(g) = 0$
otherwise.
If the difference between the largest and the second largest eigenvalues
of $T_{\psi_j}(G_j,A_j)$ is bounded from below by a constant $\epsilon>0$, 
then $\{\Gamma(G_j,A_j)\}_{j\in J}$ is a  
{\em family of expander graphs}. Now let $\{(G_j,A_j)\}_{j\in J}$ be 
the family 
of all pairs $(G,A)$ with $G = \SL_2(\mathbb{Z}/p \mathbb{Z})$,
$p$ varying over all primes, and $A$ varying over all sets of generators of 
$G$ with $d=|A\cup A^{-1}|$ fixed. The question of whether this is a family
of expander graphs may still be far from being answered.
We can prove a weaker property that has
certain consequences of its own.
\begin{cor}[of the main theorem]\label{cor:cory}
Let $p$ be a prime. Let $A$ be a set of generators of 
$G = \SL_2(\mathbb{Z}/ p \mathbb{Z})$. 
Let $\psi$ be a 
symmetric 
probability distribution on $G$ whose support contains $A$;
let $\eta = \min_{g\in A\cup A^{-1}} \psi(g)$.
Then the second largest eigenvalue
of $T_{\psi}(G,A)$ is at most
$1 - \frac{C}{\eta (\log p)^{2 c}}$, 
where $c$ and $C>0$ are absolute constants.
\end{cor}
Here $c$ is the same as in the main theorem.
\begin{proof}
Immediate from the main theorem and the standard bound for the spectral
gap in terms of $\eta$ and the diameter (see, e.g., \cite{DSC}, Cor.\ 1).
\end{proof}

From now on, assume for notational convenience that $A = A^{-1}$, and 
choose the following probability distribution on $G$:
\begin{equation}\label{eq:poof}
\psi(g) = \begin{cases} \frac{1}{2 |A|} \delta_A(g)
&\text{if $g$ is not the identity,}\\
\frac{1}{2 |A|} \delta_A(g) + \frac{1}{2} &\text{if $g$ is the identity,}
\end{cases}\end{equation}
where $\delta_A$ is the characteristic function of $A$.
For every positive integer $n$ and every $g_0\in G$, 
let $\phi_{n,g_0}$ be the probability
distribution on $G$ defined as a vector $\phi_{n,g_0} =
(T_{\psi}(G,A))^n \delta_{g_0}$, where the transition matrix
$T_{\psi}(G,A)$ is as before and $\delta_{g_0}$ is the characteristic function
of $g_0$ seen as a vector of length $|G|$. We may regard $\phi_{n,g_0}$
as the outcome of a so-called {\em lazy random walk}: start at a vertex $g_0$
of $\Gamma(G,A)$ and do the following $n$ times -- throw a coin into
the air, take a random edge out of your current vertex if it is heads, but
stay in place if it is tails.

The {\em mixing time} $\mix_{G,A}$ of the lazy random walk on
$\Gamma(G,A)$ is defined to be the smallest positive integer $n$ such that
\begin{equation}\label{eq:mix}
\sum_{g\in G}
\left|\phi_{n,g_0}(g) - \frac{1}{|G|}\right | \leq \frac{1}{2} .\end{equation}
It is clear that $\mix_{G,A}$ is independent of $g_0$.
The constant $\frac{1}{2}$ in (\ref{eq:mix}) is conventional;
if it were changed to $1/1000000$, the mixing time would
change by at most a constant factor.
\begin{cor}[of Corollary \ref{cor:cory}]
Let $p$ be a prime. Let $A$ be a set of generators of $G = \SL_2(\mathbb{Z}/
p \mathbb{Z})$. Then the mixing time $\mix_{G,A}$ is 
$O(|A| (\log p)^{2 c + 1})$, where $c$ and the implied constant are absolute.
\end{cor}
Again, the constant $c$ is as in the main theorem.
\begin{proof}
Immediate from Corollary \ref{cor:cory} via \cite{DSC}, Lemma 2.
(For $\psi$ as in (\ref{eq:poof}), the transition matrix
$T_{\psi}(G,A)$ has no negative eigenvalues; see \cite{DSC}, Lemma 1.)
\end{proof}

\begin{center}
* * *
\end{center}

By a {\em word} on the symbols $x_1,x_2,\dotsc,x_n$ we mean, as is usual,
a product of finitely many copies of
 $x_1, x_1^{-1}, x_2, x_2^{-1},\dotsc,x_n^{-1}$. A {\em trivial word}
is a product of finitely many terms of the form $g g^{-1}$, where $g$
is any word.
\begin{cor}[of the key proposition, part (\ref{it:agust})]\label{cor:aj}
Let $A$ be a set of generators of a free subgroup of
$\SL_2(\mathbb{Z})$. Let $p$ be any prime for which the reduction
$\bar{A}\subset \SL_2(\mathbb{Z}/p\mathbb{Z})$ of $A$ modulo $p$ generates
a free subgroup of 
$\SL_2(\mathbb{Z}/p\mathbb{Z})$. Then the diameter of the Cayley graph
$\Gamma(\SL_2(\mathbb{Z}/p\mathbb{Z}),\bar{A})$ 
is $O_A(\log p)$, where the implied constant
depends only on $A$.
\end{cor}
We may take, for example,
 $A$ as in (\ref{eq:knop}) or (\ref{eq:ujo}), with $p\geq 5$.
\begin{proof}
Let $g_1,g_2,\dotsc,g_n\in \SL_2(\mathbb{Z})$ be the elements of $A$.
Let $w(x_1,x_2,\dotsc,x_n)$ be a non-trivial word on $x_1,x_2,\dotsc,x_n$. Since $A$ generates
a free group, $w(g_1,g_2,\dotsc,g_n) \ne I$. Suppose that
$w(\bar{g}_1,\bar{g}_2,\dotsc,\bar{g}_n)$ equals the
identity in $\SL_2(\mathbb{Z}/ p \mathbb{Z})$, where $\bar{g}_1,
\dotsc, \bar{g}_n$ are the reductions mod $p$ of $g_1,\dotsc,g_n$.
Then at least one of the entries of $w(g_1,g_2,\dotsc,g_n)$ must have
absolute value at least $p-1$.
Yet it is clear that this is impossible if $w$ is of length
$\leq k \log p$, where $k>0$ is a constant depending only on $A$. 
(Cf.\ \cite{Ma}.) 

We thus have that any two distinct products of length at most
$\frac{k}{2} \log p$ on the symbols
$x_1,\dotsc,x_n$ must take distinct values in $\SL_2(\mathbb{Z}/p \mathbb{Z})$
for $x_1 = \bar{g}_1,\dotsc, x_n = \bar{g}_n$. We obtain that
$|\bar{A}^{\lfloor \frac{k}{2}\log p \rfloor}| \geq
n^{\lfloor \frac{k}{2} \log p\rfloor}$. 
For all $p$ larger than an absolute constant,
we have $n^{\lfloor \frac{k}{2} \log p\rfloor} \geq p^{\epsilon}$, 
where $\epsilon>0$ depends only on $k$, and hence only on $A$.
We apply part (\ref{it:agust}) of the key proposition to 
$\bar{A}^{\lfloor \frac{c}{2}\log p\rfloor}$,
and conclude
that $\diam(\Gamma(\SL_2(\mathbb{Z}/p \mathbb{Z}))) \leq C \log p$ for
some constant $C$ depending only on $A$.
\end{proof}

The following lemma seems to be folkloric. A more general statement was proved
in unpublished work by A. Shalev \cite{Lub}. Similar results have been
discovered independently by others; in particular, a generalization will appear in
a paper by Gamburd et al. \cite{Gam}. We give a proof for the sake of
completeness.
\begin{lem}\label{lem:shal}
Let $p$ be a prime. Let $G = \SL_2(\mathbb{Z}/ p \mathbb{Z})$.
Let $\mathscr{C}_p$ be the set of all pairs 
$(g,h)\in G^2$ such that $g$ and $h$ generate 
$G$. There is an absolute constant
$c>0$ such that $\Gamma(G,\{g,h\})$ has loops of length
$\leq c \log p$ for at most $o(|\mathscr{C}_p|)$ pairs $(g,h)\in \mathscr{C}_p$, where the rate of convergence to $0$ of $o(|\mathscr{C}_p|)$ is absolute.
\end{lem}
\begin{proof}
Let $w(g,h)$ be a non-trivial word. Let $f_{12}, f_{21} \in
\mathbb{Z}\lbrack x_1,x_2,\dotsc,x_n\rbrack$ be the upper-right and
lower-left entries of the matrix obtained by formally replacing all
occurrences of $g$, $h$, $g^{-1}$, $h^{-1}$ in $w(g,h)$ by the matrices
\[\left(\begin{array}{cc} x_1 & x_2\\ x_3 & x_4\end{array}\right),\;\;
\left(\begin{array}{cc} x_5 & x_6\\ x_7 & x_8\end{array}\right),\;\;
\left(\begin{array}{cc} x_4 & -x_2\\ -x_3 & x_1\end{array}\right),\;\;
\left(\begin{array}{cc} x_8 & -x_6\\ -x_7 & x_5\end{array}\right),\]
respectively. Either $f_{12}$ or $f_{21}$ is not identically equal to zero:
let $A$ be as in (\ref{eq:knop}), and denote its elements by $X$ and $Y$;
since $X$ and $Y$ generate a free subgroup
of $\SL_2(\mathbb{Z})$, at least one of the upper-right and
lower-left entries of $w(X,Y)$ or $w(Y,X)$ must be non-zero. (We cannot
have $w(X,Y) = -I = w(Y,X)$, and neither 
$w(X,Y) = I$ nor $w(Y,X)= I$ is possible.)

Assume henceforth that the length $\ell$ of $w$ is at most 
$\frac{\log (p-2)}{\log 2}$.
The coefficients of $f_{12}$ and $f_{21}$ are bounded above in absolute value
by $2^{\ell} \leq p-2$. Hence at least one of the reductions 
$\bar{f}_{12}, \bar{f}_{21}\in (\mathbb{Z}/p \mathbb{Z})\lbrack x_1, x_2,
\dotsc,x_8\rbrack$ is non-zero. Choose one of the non-zero reductions
 and call it $P$.

Since $P$ is a non-zero polynomial of degree at most $\ell$, 
there are at most $8 \ell p^7$ tuples $(x_1,\dotsc,x_8)\in (\mathbb{Z}/p
\mathbb{Z})^8$ such that $P(x_1,\dotsb,x_8) = 0$.
(While this follows immediately from the Lang-Weil estimates, it is
also quite easy to give an elementary proof. For every tuple 
$(x_2,\dotsc,x_8)\in (\mathbb{Z}/p \mathbb{Z})^7$, either there are no more
than $\ell$ values of $x_1$ with $P(x_1,\dotsc,x_8)=0$, or 
$f_{(1)}(x_2,\dotsc,x_8) = 0$, where 
$f_{(1)}$ is the leading coefficient of $f$ considered as a polynomial
on $x_1$. If $f_{(1)}(x_2,\dotsc,x_8)=0$, repeat the argument with
$f_{(1)}$ instead of $f$ and $(x_2,\dotsc,x_8)$ instead of
$(x_1,\dotsc,x_8)$.)
 Take any
$g,h\in \SL_2(\mathbb{Z}/p \mathbb{Z})$ such that $w(g,h)=I$. Then,
for all $c_1,c_2\in (\mathbb{Z}/ p \mathbb{Z})^*$, both the upper-right
and lower-right entries of $w(c_1 g, c_2 h)$ are $0$. Moreover, each
pair $c_1 g, c_2 h \in M_2(\mathbb{Z}/p \mathbb{Z})$ can arise
from at most four different pairs $g,h\in \SL_2(\mathbb{Z}/p \mathbb{Z})$.
Since every pair $c_1 g$, $c_2 h$ gives a distinct solution to
$P(x_1,\dotsc,x_8)=0$, there are at most $32 \ell p^5$ pairs $g, h
\in \SL_2(\mathbb{Z}/ p \mathbb{Z})$ such that $w(g,h) = I$.

There are at most $4^l + 4^{l-1} + \dotsb + 1 < 4^{l+1}$ distinct words
$w$ on $g$ and $h$ of length at most $l$. We conclude that,
for every $l \leq \frac{\log (p-2)}{\log 2}$, there are
fewer than $32 l 4^{l+1} p^5$ pairs $g,h\in \SL_2(\mathbb{Z}/ p \mathbb{Z})$
such that $w(g,h) = I$ for some non-trivial word $w$ of length at most $l$.
Set $l = \frac{\log p}{2 \log 4}$.
Our aim is to show that $32 l 4^{l+1} p^5 \ll p^{5.5} \log p$ is small
compared to $|\mathscr{C}_p|$; it will suffice to show that few of
the $((p^2 - 1) p)^2$ pairs
$(g,h)\in (\SL_2(\mathbb{Z}/ p \mathbb{Z}))^2$ are not in $\mathscr{C}_p$.

Every proper subgroup of $\SL_2(\mathbb{Z}/ p \mathbb{Z})$ is contained
in at least one of (a) $O(p)$ subgroups of $\SL_2(\mathbb{Z}/p)$ of
order $O(p^2)$, (b) $O(p^2)$ subgroups of order $O(p)$, or (c)
$O(p^3)$ subgroups of order $O(1)$, where the implied constants are absolute.
Tautologically, a pair of elements of a group $G$ fail to generate
$G$ if and only if they are both contained in some proper subgroup of $G$.
Hence there are at most $O(p^5)$ pairs $(g,h)\in 
(\SL_2(\mathbb{Z}/p\mathbb{Z}))^2$ not in $\mathscr{C}_p$.

We conclude that there are at most $O(|\mathscr{C}_p| (\log p)/p^{1/2})$
pairs $(g,h)\in
\mathscr{C}_p$ for which the graph $\Gamma(G,\{g,h\})$ has loops of length
$< \frac{\log p}{2 \log 4}$. (A trivial
change in the argument would give the bound
 $O_{\epsilon}(|\mathscr{C}_p| (\log p)/p^{1 - \epsilon})$
for $\epsilon>0$ arbitrary.)
\end{proof}

We can now answer in the affirmative a question of Lubotzky's 
(\cite{Lu}, Prob.\ 10.3.3).
\begin{cor}[of the key proposition, part (\ref{it:agust})]\label{cor:uj}
Let $p$ be a prime. Let $G = \SL_2(\mathbb{Z}/ p \mathbb{Z})$.
Let $\mathscr{C}_p$ be the set of all pairs 
$(g,h)\in G^2$ such that $g$ and $h$ generate 
$G$. There is an absolute constant
$C>0$ such that $\diam(\Gamma(G,\{g,h\})) \leq C \log p$ for all pairs 
$(g,h)\in \mathscr{C}_p$
outside a subset of $\mathscr{C}_p$ of cardinality $o(|\mathscr{C}_p|)$,
where the rate of convergence to $0$ of $o(|\mathscr{C}_p|)$,
is absolute.
\end{cor}
\begin{proof}
By Lemma \ref{lem:shal}, all pairs $(g,h)\in \mathscr{C}_p$ outside
a subset of $\mathscr{C}_p$ of cardinality $o(|\mathscr{C}_p|)$
yield graphs $\Gamma(G,\{g,h\})$ without loops of length $\leq c \log p$,
where $c>0$ is absolute. Let $(g,h)$ be any such pair. Then
$|\{g,h\}^{\lfloor \frac{c}{2} \log p \rfloor}| = 
|2^{\lfloor \frac{c}{2} \log p \rfloor}| \ll p^{\frac{c \log 2}{2}}$. (Cf.\ the proof of Cor.\
\ref{cor:aj}.) We apply part (\ref{it:agust}) of the key proposition to
$A = \{g,h\}^{\lfloor \frac{c}{2} \log p \rfloor}$ and are done.
\end{proof}

In Corollaries \ref{cor:aj} and \ref{cor:uj}, only the second part of the
key proposition was directly invoked. Of course, the proof of 
part (\ref{it:agust}) of the key proposition does use part
(\ref{it:wark}), but only with
 $|A| > p^{\delta}$, where $\delta>0$ is fixed. This means in turn
that the sum-product estimate (Theorem \ref{thm:bk}) is used
only for subsets of $\mathbb{F}_q^*$ whose cardinality is greater
than $p^{\epsilon}$, where $\epsilon>0$ is fixed. Thus, the results
in \cite{Ko} are not used. Since the sum-product estimates in
\cite{BKT} are purely combinatorial, the proofs of Cor.\ \ref{cor:aj}
and \ref{cor:uj} are ultimately free of arithmetic.

{\em Note added in proof.} (a) Bourgain and Gamburd have recently derived
results much stronger than Corollaries \ref{cor:aj} and \ref{cor:uj}
from the key proposition of the present paper; see \cite{BG}. (b) There
is now a proof (\cite{TV}, \S 2.8) of the sum-product theorem that does not 
involve Stepanov's
method even for subsets of $\mathbb{F}_q^*$ of cardinality smaller than
$p^{\epsilon}$. Thus, all that is not additive
combinatorics has disappeared from what is employed in this paper. 

\end{document}